\documentstyle [12pt]{article}
\newfont{\bbb} {msbm10}
\newcommand{\Bbb}[1]{\mbox{\bbb#1}}
\newcommand{\R}{\Bbb{R}}
\newcommand{\bS}{\Bbb{S}}

\newcommand{\T}{\Bbb{T}}
\newcommand{\Z}{\Bbb{Z}}
\newcommand{\CT}{{\cal{T}}}
\newcommand{\sm}{\setminus}
\newcommand{\sbs}{\subset}
\newcommand{\ra}{\rightarrow}
\begin{document}
\title{Some Remarks on the Geodesic Completeness of Compact Nonpositively 
Curved Spaces}
\author{Pedro Ontaneda}
\maketitle

Let $X$ be a geodesic space.
We say that $X$ is {\ geodesically complete} if every geodesic segment
$\beta :[0,a]\rightarrow X$ from $\beta (0)$ to $\beta (a)$ can be 
extended to a geodesic ray $\alpha : [0,\infty )\rightarrow X$, 
(i.e. $\beta (t)=\alpha (t)$, for $0\leq t\leq a$).

If $X$ is a compact npc space (``npc" means: ``non-positively curved")
then it is  {\it almost geodesically complete}, see \cite{O}. 
($X$, with metric $d$, is almost geodesically complete if its universal cover $\tilde{X}$
satisfies the following property:
there is a constant $C$
such that  for every $p,q\in \tilde X$ there is a geodesic ray $\alpha :
[0,\infty)\rightarrow \tilde X$, $\alpha (0)=p$, and $d(q,\alpha)\leq C$.)

Then it is natural to ask if in fact every compact npc
space has some kind of geodesically complete npc core.
In view of this, we ask the following question: 
(this question was already stated  in \cite{B}, p. 4)
\vspace{.2in}

{\it Question.} {\it Does a compact npc geodesic space $Z$ have a 
subspace $A$, such that
 the inclusion $A\ra Z$ is a homotopy equivalence and $A$, with the 
intrinsic metric, is npc and geodesically complete?}
\vspace{.2in}

{\it Remark.} 
One could also ask whether there is such an $A$ that is also 
totally
geodesic in $Z$, but this is certainly false. Just take $Z$ to be a 
flat 2-torus
minus an open disk. Then the only possibility for $A$ is a figure eight 
that can not
be totally geodesic in $Z$.
\vspace{.2in}

If $Z$  is piecewise smooth npc 2-complex, then, by collapsing the
free faces, we get a subcomplex which is, with its intrinsic metric,
geodesically complete and npc (see  \cite{BB}, and also \cite{BH}, p.208).
Hence the answer to the
problem above is YES, when $Z$ is piecewise smooth  a 2-complex.

For the general case the answer is NO. To show this,
in section 2 we construct a space $Z$ with the following properties:
\vspace{.2in}

{\bf First Example.} {$Z$ is a finite piecewise flat npc 3-complex 
with the property that there is no $PL$ subspace $A$ of $Z$ such that

(i) The inclusion $\iota : A\ra Z$ is a homotopy equivalence.

(ii) $A$, with the intrinsic metric determined by the metric of
$Z$, is npc geodesically complete.}
\vspace{.2in}

But the space $Z$ we construct in this example {\it does have} a 
subspace $A$
satisfying $(i)$ and also $(ii)$ but with other geodesic metric (not 
the intrinsic one).
Hence we are forced to ask a deeper question.
\vspace{.2in}

{\it Question.} {\it Does a compact npc geodesic space $Z$ have a 
subspace $A$, such that
 the inclusion $A\ra Z$ is a homotopy equivalence and $A$ admits a 
geodesic metric that  is npc and geodesically complete?}
\vspace{.2in}

For this question the answer is also NO. To show this, in section 3 
we construct an space $Z$ with the following properties:
\vspace{.2in}

{\bf Second Example.} {\it $Z$ is a finite piecewise flat npc 3-complex 
with the property that there is no $PL$ subspace $A$ of $Z$ such that

(i) The inclusion $\iota : A\ra Z$ is a homotopy equivalence.

(ii) $A$ admits a npc geodesically complete geodesic metric.}
\vspace{.3in}

Of course the second example is also acounterexample to the
first question but,
eventhough it is a little non-standard to present these counterexamples
in this way, we think it is interesting  to show how the constructions evolve.

Still, the space $Z$ we construct for theorem E {\it is} homotopically 
equivalent 
to a npc geodesically complete space. Hence we are again forced to ask 
an even deeper
question, for which we do not know the answer.
\vspace{.2in}

{\it Open Question.} {\it Is a compact npc geodesic space $Z$ 
homotopically equivalent to a compact npc geodesically complete 
geodesic space $Z'$?}
\vspace{.2in}

We say that a group is npc if it is isomorphic to the fundamental group 
of a
compact npc geodesic space, and a group is gc-npc (i.e. geodesically 
complete
non-positively curved) if it is isomorphic to the fundamental group of 
a  compact npc geodesically complete geodesic space. With these 
definitions, we restate the open problem above in the following way: Is every npc 
group a gc-npc group? Remark that if we drop the compactness condition the 
answer is always YES, at least in the $PL$ category: let $K$ be a piecewise flat npc 
simplicial complex. Embed it in some $\R^N$ and let $Y$ be the hyperbolization
of $\R^N$ relative to $K$. $Y$ can be given a complete npc metric 
having $K$ as a totally geodesic subspace (see \cite{BHu}). Let $Z$ be the cover 
of $Y$ relative to the subgroup $\pi_{1}(K)$ of $\pi_{1}(Y)$. Then $Z$ is a non-compact
npc geodesically complete space having $K$ as a deformation retract.

For a group $\Gamma$ define the {\it min-gc}-dimension and the {\it 
max-gc}-dimension in the following way. {\it min-gc-dim} $ \Gamma$ is the minimun of the 
$n$'s such that there is a compact npc geodesically complete $n$-complex with 
fundamental group isomorphic to $\Gamma$. If there is not such complex we put
{\it min-gc-dim} $ \Gamma = +\infty$. Similarly  {\it max-gc-dim} $ 
\Gamma$ is the maximun of the $n$'s such that there is a compact npc geodesically complete 
$n$-complex with fundamental group isomorphic to $\Gamma$. If there is not such complex 
we put {\it max-gc-dim} $\Gamma = -\infty$, and if there are infinite such 
$n$'s we write {\it max-gc-dim} $\Gamma = +\infty$

Note that we always have $gd\,\Gamma\leq $ {\it min-gc-dim} $\Gamma$, 
provided {\it min-gc-dim} $\Gamma\neq +\infty$. Here $gd\,\Gamma$ is the 
geometric  dimension of $\Gamma$.
This motivates the following quetion.
\vspace{.2in}

{\it Open Question.} {\it Is it true that min-gc-dim $\Gamma$ = gd $\Gamma$, provided that 
min-gc-dim $\Gamma\neq +\infty$?}
\vspace{.2in}

In \cite{Br2} it is given an example of a group $\Gamma$, with $gd\, 
\Gamma=2$, that is the fundamental group of a npc compact 3-complex but 
it is not the fundamental group of a npc 
compact 2-complex. Hence, at least one of the two open problems above 
have a negative answer.

It is easy to show that for the infinite cyclic group $F_{1}$ we have
{\it min-gc-dim} $F_{1}$={\it max-gc-dim} $F_{1}=gd \, F_{1}=cd =\, 
F_{1}=1$. Here $cd$ is the
cohomological dimension.
But for $n>2$ things are different. ($F_{n}$ denotes the free group 
with $N$ generators.)
\vspace{.2in}

{\bf Proposition 1.} {\it max-gc-dim $F_{n}>1$, for $n\geq 2$.}
\vspace{.2in}

Because {\it min-gc-dim $F_{n}$=1}, for all $n$, we have that in 
general is not true that
{\it min-gc-dim} and {\it max-gc-dim} are equal. And also, in general, 
it is not
true that {\it max-gc-dim} and {\it gd} are equal. 
In fact in general it is not true that
{\it min-gc-dim} and {\it max-gc-dim} are equal, even if $\Gamma$ is 
the fundamental group of a closed manifold, as the next proposition 
asserts.
\vspace{.2in}

{\bf Proposition 2.} {\it Let $G_{n}$ denote  the fundamental group of 
the genus n surface.
Then max-gc-dim $G_{n}>2$, for n large.}
\vspace{.2in} 

In fact this proposition is a corollary of the following proposition.
\vspace{.2in}

{\bf Proposition 3.} {\it For every m there is an k such that 
max-gc-dim $G_{n}\geq m$,
for $n\geq k$.}
\vspace{.2in}

In fact, using the same method of the proof of Proposition 3
we get that also for $F_{n}$ we get:
For every m there is an k such that max-gc-dim $F_{n}\geq m$,
for $n\geq k$.

From the proof of  proposition 1, it seems that the answer to the 
following problem is YES.
\vspace{.2in}

{\it Open Question.} {\it Is max-gc-dim $F_{n}$ finite?}
\vspace{.2in}

Or, more generally.
\vspace{.2in}

{\it Open Question.} {\it Is max-gc-dim $\Gamma$ finite, provided 
min-gc-dim $\Gamma\neq +\infty$?}
\vspace{.2in}

{\it Open Question.} {\it Assume that min-gc-dim $\Gamma\neq +\infty$, 
and let
min-gc-dim $\Gamma \leq$ m $\leq$ max-gc-dim $\Gamma$. Does there exist 
a compact npc geodesically complete m-complex with fundamental
group isomorphic to $\Gamma$?}
\vspace{.2in}

The next proposition says that a npc group can be made geodesically
complete after making a free product with some free group, at least in 
the
$PL$ category. 
\vspace{.2in}

{\bf Proposition 4.} {\it Let $\Gamma$ be the fundamental group of
a finite piecewise flat npc simplicial complex. Then there is an n
such that $\Gamma * F_{n}$ is gc-npc.}
\vspace{.2in}

The number $n$ in proposition D is very large and it would be 
interesting
to know the relationship between $n$ and $\Gamma$. For example for 
$n=1$
we can ask:
\vspace{.2in}

{\it Open Question.} {\it Does $\Gamma * F_{1}$  gc-npc implies $\Gamma$  
gc-npc.}
\vspace{.3in}

Finally, it is important to remark that geodesic completeness
is a useful property related to rigidity results. In certain cases 
rigidity
results in the manifold category can be generalized replacing the 
manifold
condition by the geodesic completeness condition. See for example the
work of Leeb \cite{L}, and Davis-Okum-Zheng \cite{DOZ}.
One of the simplest cases of this is 
the topological version of Gromoll-Wolf-Lawson-Yau torus theorem 
(see \cite{B},\cite{Br1}). For example, if $Z$ is a compact npc geodesic 
space
homotopically equivalent to some torus $\T^n$ then there is a totally
geodesic embedding $\T^{n}\ra Z$ and easily follows that if $Z$ is 
geodesically
complete, $Z$ is isometric to $\T^n$. But propositions C and D show 
that geodesic
completeness is, in general, not enough. In fact, the dimension of the 
spaces
could be the same and, still, we do not get rigidity, as the next 
proposition
shows.
\vspace{.2in}

{\bf Proposition 5.} {\it For $n\geq 2$ there is a finite piecewise 
flat npc geodesically
complete 2-complex $Z$ with $\pi_{1}(Z)\cong G_{n}$ and $Z$ is not a 
2-manifold.}
\vspace{.2in}

Hence, if we want to generalize rigidity results (like Farrell-Jones 
topological
rigidity of npc manifolds) to non-manifold categories, we need more 
than geodesic
completeness. On the other hand the singularities of the $Z$ we 
construct in proposition
F are quite trivial and one wonders whether this is a general fact.
\vspace{.3in}

Here is a short outline of the paper. In section 1 we give some 
definitions and a lemma. In section 2 we construct the first example.
In section 3 we construct the second example. In section 
4 we prove propositions 1,3,4 and 5.
\vspace{.3in}

{\bf 1. Definitions.}
\vspace{.2in}

{\bf 1. Definitions and Lemmas.}
\vspace{.2in}

Let $\Delta^n$ denote the canonical $n$-simplex (i.e. $\Delta^n$  is 
the convex hull
of the $n+1$ points $(1,0,...,0),...,(0,...,0,1)$ in $\R^{n+1}$).
Let also $\T^n$ denote the $n$-torus with its canonical $PL$ structure.
Now, let $\Lambda$ be a $PL$ subspace of $\T^3$ $PL$ homeomorphic to
$\Delta^2$. We construct the $PL$ space $T$ by taking two copies $T^0$ 
and $T^1$ of
$\T^3$ and identifying them along $\Lambda$, and we consider $\Lambda$, 
$T^0$ and $T^1$
as a $PL$ subspaces of this quotient space $T$.

Consider a pair $(K,A)$, where $K$ is a $PL$ space, 
$A=\{\sigma_{1},...,\sigma_{r}\}$,
and each $\sigma_i$ is a $PL$ subspace of $K$ $PL$ homeomorphic to 
$\Delta^2$.
Given these data we construct the $PL$ space $\CT (K,A)$ by taking, for 
each $\sigma_i$,
a copy of $T$ and identifying $\Lambda$ with $\sigma_i$.
We consider $K$, each copy $T_{i}$ of $T$,
and each copy $T_{i}^{j}$ of $\T^j$, $j=0,1$ as being $PL$ subspaces
of $\CT (K,A)$.
Also, if $J$ is a $PL$ subspace of $K$ with $\sigma_{i}\sbs J$, for all 
$i$,
we can consider $\CT (J,A)\sbs \CT (K,A)$.

Finally, a free face of a simplicial complex $J$ is a $n$-simplex 
$\sigma$ in $J$
that is the face of exactly one $n+1$-simplex in $J$.
\vspace{.2in}

{\bf Lemma 1.1.} {\it Let $K$ and $A$ be as above and $Y$ a closed $PL$ 
subspace
of $X =\CT (K,A)$ such that the inclusion
$\iota : Y \rightarrow X$ is a homotopy equivalence.
Then  $T_{i}\sbs Y,
\,\, \, i=1,...,r.$}
\vspace{.2in}

{\bf proof.} Take a point $p$ in $T_{i}\setminus K$, for some 
$i=1,...r$.
We have that $H_{3}(X , X\sm\{ p\} )\cong \Z$. (Here $H$ denotes
homology with $\Z$ coefficients.)
Also, the map $H_{3}(X )\rightarrow H_{3}(X , X\sm\{ p\} )$ is onto.
Because $\iota_{*}:H_{3}(Y)\ra H_{3}(X )$ is an isomorphism,
we get that the composition $H_{3}(Y)\ra
H_{3}(X )\rightarrow H_{3}(X , X\sm\{ p\} )$ is onto.
Hence $Y$ is not contained in $X\sm \{ p\}$.
This means that $T_{i}\sm K\sbs Y$. Hence $T_{i}={\overline{T_{i}\sm 
K}}\sbs\overline{Y}=Y$.
This proves the lemma.
\vspace{.2in}

{\bf Remark.} Note that lemma 5.1 implies that each $\sigma_i$ is in 
$Y$.
Hence we have that $Y=\CT (J,A)$, for some closed $PL$ subspace $J$ of 
$K$,
with $\sigma_{i}\sbs J$, $i=1,...,r$.
\vspace{.2in}

{\bf Lemma 1.2.} {\it Let $K$ and $A$ be as above and $J$ a $PL$ 
subspace
of $K$. If the inclusion $\eta :\CT (J,A)\ra \CT (K,A)$ is a homotopy
equivalence, then the inclusion $\iota :J\ra K$ is also a homotopy 
equivalence.}
\vspace{.2in}

{\bf proof.} We have the following diagram.

$$\begin{array}{ccc}

J &\stackrel{\iota}{\longrightarrow} & K \\
 \alpha\downarrow\,\,\,\,&&\beta\downarrow\,\,\,\, \\
\CT (J,A)&\stackrel{\eta}{\longrightarrow} &\CT (K,A)\\
r\downarrow\,\,\,\,&& s\downarrow\,\,\,\, \\
J &\stackrel{\iota}{\longrightarrow} & K
\end{array}$$

Here all the arrows of the upper square are inclusions and the vertical 
arrows
of the lower square are retractions induced by retractions of each 
$T_i$ to
$\sigma_{i}$. Applying the $i$-th homotopy functor to the diagram above
we get

$$\begin{array}{ccc}

\pi_{i}(J) &\stackrel{\iota_{*}}{\longrightarrow} &\pi_{i}( K) \\
 \alpha_{*}\downarrow\,\,\,\,&&\beta_{*}\downarrow\,\,\,\, \\
\pi_{i}(\CT (J,A))&\stackrel{\eta_{*}}{\longrightarrow} &\pi_{i}(\CT 
(K,A))\\
r_{*}\downarrow\,\,\,\,&& s_{*}\downarrow\,\,\,\, \\
\pi_{i}(J )&\stackrel{\iota_{*}}{\longrightarrow} &\pi_{i}( K)
\end{array}$$

Because $r\alpha =id_{J}$, $s\beta =id_{K}$ and $\eta_{*}$ is, by 
hypothesis,
an isomorphism we get that $\iota_{*}$ is an isomorphism. This proves 
the lemma.
\vspace{.2in}

{\bf Lemma 1.3.} {\it Let $J$ be a $PL$ subspace of $\Delta^3$,
$A=\{\sigma_{1},...,\sigma_{r}\}$ a set of $PL$ subsets of $J$, $PL$
equivalent to $\Delta^2$. Assume that $B=\cup A=\cup_{i}\sigma_{i}$
is contractible. If $Y=\CT (J,A)$ admits a npc geodesically complete 
geodesic
metric, then $J$ is two dimensional
(i.e. $PL$ equivalent to a two complex).}
\vspace{.2in}

{\bf proof.} We can assume that $J$ is a simplicial complex such that
each $\sigma_i$ is a subcomplex of $J$.
Let $K$ be a simplicial complex $PL$ equivalent to $\Delta^3$,
such that we can consider $J$ as being a subcomplex of $K$.
Suppose that $J$ contains at least one 3-simplex.
Let $c=\sum c_{k}$, a $\Z_{2}$ simplicial 3-chain in $\Delta^3$, where 
the $c_i$'s
are all the 3-simplices of $J$, without repetition. Because $Y$ is
geodesically complete we have that $Y$ does not have free faces
(see \cite{O}). Hence, every 2-simplex in $J$ is in the boundary
of at least two 3-simplices of $Y$. But $J\sbs\Delta^3$, thus
every 2-simplex of $J$ is in the boundary of exactly two
3-simplices, unless the 2-simplex lies in some $\sigma_{i}\sbs B$.
Hence $\partial c =0$, and follows that $H_{3}(\Delta^{3},B)\neq 0$.
This is a contradiction because $B$ is contractible, hence
$H_{3}(\Delta^{3},B)=0$. This proves the lemma.
\vspace{.2in}

{\bf Lemma 1.4.} {\it Assume that $X=\CT (K,A)$ admits a npc
geodesic metric $d$. Then $K$ is totally geodesic in $X$.}
\vspace{.2in}

{\bf proof.} First, by the topological flat torus theorem (see, for 
example, \cite{Br1}), it is easy to see
that each $T_{i}^j$, is totally geodesic, for all $i$ and $j=0,1$.
Hence $\Lambda_{i}=T_{0}^{i}\cap T_{i}^1$ is totally geodesic.
We prove now that $K$ is also totally geodesic.
Let $p,q\in K$, and $\alpha : [0,d(p,q)]\ra X$ a geodesic with
$\alpha (0)=p$, and $\alpha (d(p,q))=q$. For each $T_i$ let
$a=min\{ t\,  :\,\alpha (t)\in T_{i}\}$ and $b=max\{ t\, :\,\alpha 
(t)\in T_{i}\}$.
Then $\alpha\mid_{[a,b]}$ is a geodesic joining $\alpha (a)$ to $\alpha 
(b)$.
But $\alpha (a), \alpha (b)\notin T_{i}\sm K$ (because is open), hence
$\alpha (a), \alpha (b)\in T_{i}\cap K=\Lambda_i$. But we saw before 
that
$\Lambda_i$ is totally geodesic. Consequently $\alpha [a,b]\sbs 
\Lambda_{i}\sbs K$.
This proves the lemma.
\vspace{.4in}

{\bf 2. First Example.}
\vspace{.2in}

Consider $\Delta^3\sbs\R^4$ with the flat metric induced by
$\R^4$. Fix a vertex $v$ of $\Delta^3$. Let $A=\{
\sigma_{1},\sigma_{2},\sigma_{3}\}$, where the $\sigma_i$'s are
the 2-simplices of $\partial \Delta^3$ that contain $v$. Write
$B=\cup A= \sigma_{1}\cup\sigma_{2}\cup\sigma_{3}$. Take $X=\CT
(\Delta^{3}, A)$. Consider $\T^3$ with it canonical flat metric
and $\Lambda\sbs\T^3$ isometric to $\Delta^2$. Because each
$\sigma_i$ in $\Delta^3$ is isometric to $\Delta^2$, all the
identifications used to construct $X$ are isometries. Hence this
determines a piecewise flat npc metric on $X$.
\vspace{.2in}

{\bf Proposition 2.1.} {\it There is no $PL$ subspace $Y$ of $X$ such 
that

(i) The inclusion $\iota : Y\ra X$ is a homotopy equivalence.

(ii) $Y$, with the intrinsic metric determined by the metric of
$X$, is geodesically complete non-positively curved.}
\vspace{.2in}

{\bf proof.} Suppose such a $Y$ exists. By the remark after lemma
5.1. we have that $Y=\CT (J,A)$, for some $PL$ subspace $J$ of
$\Delta^3$, and each $\sigma_{i} \sbs J$. By lemma 5.3,  $J$
is  two complex. Hence $link(v,J)$ is one dimensional. Because
$\sigma_{i}\sbs J$ and $v\in \sigma_i$, $i=1,2,3$, we get a loop
of length $3\frac{\pi}{3}=\pi$ in $link(v,J)$ (each $\sigma_i$
determines a path of length $\frac{\pi}{3}$ in $link(v,J)$). Thus
$J$, with the intrinsic metric, is not npc. But
lemma 5.4 says that $J$ is totally geodesic in $Y$, hence $J$ is
npc, a contradiction. This proves the
proposition. \vspace{.2in}

If we choose $J=\sigma_{1}\cup\sigma_{2}\cup\sigma_{3}$, and
$Y=\CT (J,A)$, then we can provide this $Y$ with a geodesically
complete piecewise flat npc metric: give each
$\sigma_i$ a flat metric in such a way that the angle at
$v\in\sigma_i$ is $\frac{\pi}{3}$ and choose the 2-simplices
$\Lambda_i$ accordingly. (Note that this metric on $Y$ {\it is not
the intrinsic metric}.)With this choice of $Y$ we certainly have
that the inclusion $\iota :Y\ra X$ is a homotopy equivalence.
Consequently, our space $X$ does admit a $PL$ subspace $Y$ that
admits a geodesically complete non-positively curved metric. In
our next example this will not happen. \vspace{.4in}

{\bf 3. Second Example.}
\vspace{.2in}

Consider again $\Delta^3\sbs\R^4$ with the flat metric induced by
$\R^4$. Let $L\sbs\Delta^3$ be a contractible two complex with no
free faces. (Take for example the house with two rooms, see \cite{C}.)
Let $A$ be the set of all 2-simplices of $L$ and take $X=\CT
(\Delta^{3},A)$. Note that we can give $X$ a npc
piecewise flat metric in the following way. Let $\sigma_{i}\in A$
and choose $\Lambda_{i}\in T_i$ isometric to $\sigma_i$. In this
way all the identifications used to construct $X$ are isometries
and this determines a piecewise flat npc metric
on $X$.
\vspace{.2in}

{\bf Proposition 3.1.} {\it There is no $PL$ subspace $Y$ of $X$
such that

(i) The inclusion $\iota : Y\ra X$ is a homotopy equivalence.

(ii) $Y$ admits a geodesically complete non-positively curved
geodesic metric.} \vspace{.2in}

{\bf proof.} Suppose such a $Y$ exists. By the remark after lemma
5.1. we have that $Y=\CT (J,A)$, for some $PL$ subspace $J$ of
$\Delta^3$, and that each $\sigma_{i} \sbs J$. Hence we can assume
that $J$ is a simplicial complex and that $L$ is a subcomplex of
$J$. By lemmas 5.2 and 5.3, $J$ is a contractible two complex.
\vspace{.1in}

{\bf claim.} {\it $J$ has no free faces (or, in this case, free
edges).}
\vspace{.1in}

If $J$ has a free edge $e$, we have two possibilities. First, if
$e$ is not in $L$ then $e$ would also be a free edge of $Y$, but
this is impossible because $Y$ is geodesically complete
non-positively curved, thus has no free faces (see \cite{O}). Second,
$e$ can not be in $L$ because $L\sbs J$ and $L$, by hypothesis,
has no free faces. This proves the claim. \vspace{.1in}

But the claim and lemma 5.4 imply that $J$ is a contractible
complex with no free faces that admits a geodesically complete
npc geodesic metric. This is impossible (see
\cite{O}). This proves the proposition. \vspace{.2in}

So, our $X$ above does not have a subspace that admits a
geodesically complete npc geodesic metric. But
there is a space $Z$ homotopically equivalent to $Y$ that admits a
geodesically complete npc geodesic metric. In
fact, $X$ is homotopically equivalent to a finite wedge of tori,
which certainly admits such a metric.
\vspace{.8in}

{\bf 4. Proofs of Propositions 1,3,4 and 5.}
\vspace{.2in}

{\bf proof of proposition 1.} Let $C=[0,1]\times\bS ^1$. Then $C$ is 
npc
and $S_{0}=\{ 0\}\times\bS ^{1}$ and $S_{1}=\{ 1\}\times\bS ^{1}$ are 
totally geodesic
in $C$. Let $s_{0}$ be an embedded segment of length one in the 
interior
of $C$. Let $n\geq $3. Let $s_{1},...,s_{n-2}$ be embedded segments in 
the
interior of $C$ of length $\frac{1}{n-2}$ such that 
$s_{0},s_{1},...,s_{n-2}$
are all disjoint. Let $u_{1},...,u_{n-2}$ be segments of length 
$\frac{1}{n-2}$
embedded in $S_{1}$ with $S_{1}=\bigcup u_{i}$. Because all $s_{i}$'s, 
$i>0$,
and $u_{i}$'s are totally geodesic and  isometric (they have the same
length) we can identify each $s_{i}$, $i>0,$ with $u_{i}$. Also, 
because
$S_{0}$ and $s_{0}$ have the same length there is a  surjective local
isometry $s_{0}\ra S_{0}$, and we can use this local isometry to 
identify
$s_{0}$ with $S_{0}$. Let $Z$ be the resulting space. Then $Z$ is npc,
compact and it is easy to see that $\pi_{1}(Z)\cong F_{n}$. Also the
simplicial complex $Z$ has no free faces. Hence it is geodesically
complete (see \cite{BH}, p.208).

For $n=2$ the construction is similar. Just replace the cylinder $C=[0,1]\times\bS ^1$
by a Moebius band.
This proves proposition 1.
\vspace{.4in}

{\bf proof of proposition 4.} 
The proof is similar to the proof of proposition 1. Let $Z$ be a finite 
piecewise flat npc
simplicial complex. Let $\sigma$ be a free face of dimension $l$. Then 
$\sigma$ is in the boundary of exactly one $(l+1)$-simplex. Call this 
simplex $\tau$. Subdivide $\sigma$ and
let $s_{1},...,s_{k}$ be the $l$-simplices of this subdivision and 
assume
$diam\, s_{i}<diam\, \sigma$. Then we can find disjoint  $l$-simplices 
$u_{1},...,u_{k}$ in the interior of $\tau$ such that each $u_{i}$ is 
isometric
to $s_{i}$. Identify each $s_{i}$ with $u_{i}$. 
After this there are no free faces contained in $\sigma$. Do this for 
every free face.
After a finite number of steps we get a finite piecewise flat npc 
simplicial
complex $Y$ with no free faces. Consequently $Y$ is geodesically 
complete.
Also it is easy to see that $\pi_{1}(Y)=\pi_{1}(Z)*F_{N}$ for some 
large $N$.
This proves proposition E.
\vspace{.4in}

{\bf proof of proposition 5.}
Let $n\geq$ 2. Let $P\subset\R^2$ be a regular $4n$-gon.
We can identify the sides of $P$ to obtain 
a piecewise flat npc 2-simplex $Y$ homeomorphic to the surface of genus 
$n$.
In this particular case the npc metric induced on $Y$ by $P$ is a flat 
metric,
except for the vertex $v$, where the length $\ell$ of link is larger 
than $2\pi$.
By the $PL$ Gauss-Bonnet formula $2\pi (2-2n)=2\pi -\ell$. Thus $\ell 
=2\pi (2n-1)\geq 6\pi$.
Let $s_{0}$ and $s_{1}$ be embedded segments on $Y$ of the same length 
such that 
$s_{0}\cap s_{1}=v$ and the two angles determined by $s_{0}$ and 
$s_{1}$ at $v$ are
each larger than $2\pi$. This is possible because $\ell\geq 6\pi$.
Let $Z$ be obtained from $Y$ identifying $s_{0}$ with $s_{1}$.
This space $Z$ satisfies the statement of the proposition. This proves
proposition F.
\vspace{.4in}

Before proving proposition 3 we prove a lemma.
\vspace{.2in}

{\bf Lemma 4.1.} {\it Let $a_{1},...,a_{j},\theta_{1},...,\theta_{j}$ 
be
positive real numbers. Then there is piecewise flat npc complex $W$ 
homeomorphic
to a 2-surface with an open disk deleted, such that the boundary of $W$ 
consists
of $j$ geodesic segments $s_{1},...,s_{j}$ of lengths $a_{1},...,a_{j}$ 
and the angle
at the initial point of $s_{i}$ is larger than $\theta_{i}$, 
$i=1,...j$.}
\vspace{.2in}

{\bf proof.} 
Let $Y$ be the space of the proof of proposition F. Let $V$ be obtained 
from $Y$
in the following way. Delete the segment $s_{0}$ from $Y$ and glue back 
two copies
of $s_{0}$ which will intersect in their initial and final points. Then 
$V$ is
homeomorphic to a surface with an open disk deleted. Also $V$ is 
piecewise flat
and npc and its boundary consists of two segments and two vertices with 
angles
$2\pi$ and $\ell =2\pi (2n-1)$ (see the proof of prop. F). 
Note that we can choose $Y$, $n$, and the length of $s_{0}$ 
arbitrarily. Hence taking
several spaces like this $V$ and gluing them along the boundaries it is 
not
difficult to show that we can find the required $W$.
\vspace{.2in}

{\bf proof of Proposition 3.}
Let $m\geq 3$ be a positive integer. Let $\Delta^m$ be the canonical 
$m$-simplex.
By proposition E there is a finite piecewise flat geodesically
complete npc $m$-complex $U$ with 
$\pi_{1}(U)=\pi_{1}(\Delta^{m})*F_{r}=F_{r}$, for
some large $r$, and we can assume $r$ even. Let $R$ be a 1-complex 
$PL$-embedded
in $U$ such that the inclusion $R\ra U$ is a homotopy equivalence. Then 
$R$ has the homotopy type of a wedge of $r$ circles. Because $r$ is 
even there is a map $\rho$ from the boundary
$\bS^{1}=\partial D$ of the 2-disk $D$ to $R$ such that if we glue $D$ 
with the
mapping cylinder of $\rho$ along $\partial D$ we obtain a closed 
surface of genus
$\frac{r}{2}$. We can assume that $\rho$ is a simplicial map and let
$s_{1},..., s_{j}$ be a subdivision of $\partial D$ in segments such 
that $\rho$ 
is simplicial on each $s_{i}$. Let $a_{i}$ be the length of $\rho 
(s_{i})$ and 
choose $\theta_{i}\geq 2\pi -\alpha_{i}$,
where $\alpha_{i}$ is the angle at the initial point of $\rho (s_{i})$
between the segments $\rho (s_{i})$ and $\rho (s_{i-1})$. Let $W$ be 
the space
given by lemma 8.1 corresponding to these $a_{i}$'s and $\theta_{i}$'s.
Hence we have a map $\psi :\partial W\ra U$ with $\psi (\partial 
W)\subset R$,
that is a local isometry. Let $Z$ be obtained by gluing $W$ and the 
mapping cylinder
of $\psi :\partial W\ra U$ along $\partial W$. Then $Z$ is an 
$m$-complex
and it is homotopically equivalent to  a closed surface, and it is also
a finite geodesically complete npc piecewise flat complex. This proves 
proposition 3.


\begin{thebibliography}{99}

\bibitem{B} W. Ballmann, {\em Lectures on spaces of non-positive 
curvature}, DMV Seminar,
Band 25, (1995) Birkhauser.
\bibitem{BB} W. Ballmann and M. Brin, {\em Polygonal complexes and 
combinatorial group
theory}, Geom. Dedicata 50 (1994), 165-191.


\bibitem{Br1} M. Bridson, {\em On the existence of flat planes in 
spaces of non-positive curvature},Proc. Am. Math. Soc. 123, 223-235  (1995).
\bibitem{Br2} M. Bridson, {\em Length functions, curvature and the 
dimension of
discrete groups}, preprint.





\bibitem{BH} M. Bridson and A. Haeflinger, {\em Metric spaces of 
non-positive curvature},
Springer-Verlag (1999).




\bibitem{C} M. M. Cohen, {\em a course in simple homotopy theory}, 
Springer Verlag, New York,
1972.




\bibitem{DOZ} M.Davis, B.Okum and F. Zheng, {\em Piecewise euclidean 
structures and Eberlein's
rigidity theorem in the singular case}, Geom. Topol. 3, 303-330 (1999).


\bibitem{BHu} B. Hu, {\em Whitehead groups of nonpositively curved 
polyhedra}, J. Differential Geom.
38  (1993)  501-517.

\bibitem{L} B. Leeb, {\em A characterization of irreducible symmetric 
spaces and euclidean
buildings of higher rank by their asymptotic geometry}, preprint.
\bibitem{O} P. Ontaneda, {\em Cocompact CAT(0) spaces are almost
geodesically complete}, preprint.
\bibitem{RS} C. P. Rourke and  B. J. Sanderson, {\em Introduction to
piecewise-linear  topology}, Springer-Verlag,
(1972).

\end{thebibliography}
\end{document}